\documentclass[runningheads]{lncse}

\usepackage{amsmath}
\usepackage{amsfonts}
\usepackage{epsfig}
\usepackage{graphics} 
\usepackage{graphics}
\usepackage{color}
\usepackage{manfnt}
\usepackage{mathtools}
\usepackage{multirow}
\usepackage{array}

\graphicspath{{./figures/}}

\newcommand{\A}{\mathcal{A}}
\newcommand{\B}{\mathcal{B}}

\newcommand{\micron}{\mu\text{m}}
\newcommand{\PI}{\Pi_R}

\newcommand*{\Scale}[2][4]{\scalebox{#1}{$#2$}}%
\newcommand{\reals}{\mathbb{R}}

\newcommand{\half}{\Scale[0.5]{\tfrac{1}{2}}}
\newcommand{\nhalf}{\Scale[0.5]{-\tfrac{1}{2}}}

\begin{document}

\title{Sub-voxel perfusion modeling in terms of coupled 3$d$-1$d$ problem}

\titlerunning{Perfusion modeling by coupled 3$d$-1$d$ system}

\author{ Karl Erik Holter\inst{1, 3} \and Miroslav Kuchta\inst{2} \and Kent-Andr{\' e} Mardal\inst{2, 3}}

\authorrunning{Holter, Kuchta, Mardal}   

\institute{
  University of Oslo, Department of Informatics {\tt holter@simula.no}
  \and
  University of Oslo, Department of Mathematics, {\tt mirok@math.uio.no}
  \and
  Simula Research Laboratory, \tt{kent-and@simula.no}
}

\maketitle

\begin{abstract}
  We study perfusion by a multiscale model coupling diffusion in the
  tissue and diffusion along the one-dimensional segments representing
  the vasculature. We propose a block-diagonal preconditioner for the
  model equations and demonstrate its robustness by numerical experiments.
  We compare our model to a macroscale model by [P. Tofts, Modelling in DCE MRI, 2012].
\end{abstract}

\section{Introduction}

The micro-circulation is altered in diseases such as cancer and Alzheimer's disease, as demonstrated 
with modern perfusion MRI. In cancer, the so-called enhanced permeability and retention (EPR) effect  
describes the fact that the smaller vessels in a tumor are leaky, highly permeable vessels that 
enable the tumor cells to grow quicker than normal cells.

In Alzheimer's disease (AD), the opposite is alleged to happen. According to \cite{jack2009cerebrovascular}, hypoperfusion is a precursor to AD, and the cause of the pathological cell-level changes occuring in AD. This could also explain why various kinds of heart disease are risk factors for AD, as changes in the blood pressure would affect the perfusion of the brain.


The vasculature (e.g. idealized as a system of pipes) and the surrounding
tissue are clearly three-dimensional. However, the fact that in many applications
the radii of the vessels are negligible compared to their lengths, permits reducing
their governing equations to models prescribed on (one dimensional) curves where
the physical radius $R$ enters as a parameter. In this paper, we consider a
coupled 3$d$-1$d$ system
\begin{align}
  \frac {\partial u} {\partial t} &=  D_\Omega  \Delta_{\Omega} u + \delta_{\Gamma} \beta \lambda                  &\quad\text{ in } \Omega, \nonumber \\
  \frac {\partial \hat{u}} {\partial t} &=  D_\Gamma  \Delta_{\Gamma} \hat{u} -\beta \lambda    &\quad\text{ in } \Gamma, \label{kam_eq:strongform} \\ 
\lambda &= \beta (\PI u - \hat{u})   &\text{ in } \Gamma. \nonumber
\end{align}

Here, $\delta_{\Gamma}$ is the Dirac measure on $\Gamma$;  $u$, $\hat{u}$ are
the concentrations in the tissue domain
$\Omega$ and the one-dimensional vasculature representation $\Gamma$; and
$D_{\Omega}$, $D_{\Gamma}$ are conductivities on the respective domains. The last
equation is then a generalized Starling's law; relating $\hat{u}$ and the value
$\PI u$ of $u$ averaged over the idealized cylindrical vessel surface centered around
$\Gamma$\footnote{ Let $x \in\Gamma$ and $C_R(x)$ be a circular crossection of the vessel
  surface with a plane $\left\{y\in\mathbb{R}^3, (y-x)\cdot\tfrac{\mathrm{d}\Gamma}{\mathrm{d} s}(x) = 0\right\}$
  defined by the tangent vector of $\Gamma$ at $x$. The surface avarage $\PI u$ of $u$
  is then defined by
\[
  (\PI u)(x) = \lvert C_R(x) \rvert ^{-1}\int_{C_R(x)} u(y)\,\mathrm{d}y.
\]
}. The equation thus represents a coupling between the domains.

Variants of the system \eqref{kam_eq:strongform} have been used to study
coupling between tissue and vasculature flow in numerous applications.
In \cite{goldman2000computational}, a steady state limit of the system is
used to numerically investigate oxygen supply of the skeletal muscles.
Finite differences were used for the discretization.
Using the method of Green's functions, \cite{secomb2004green}, \cite{secomb2000theoretical}
and \cite{chapman2008multiscale} studied oxygen transport in the brain and tumors.
Transport of oxygen inside the brain was also investigated by \cite{linninger2013cerebral}
using the finite volume method (FVM) and \cite{fang2008oxygen} using the
finite element method (FEM) for the 3$d$ diffusion problem and FVM elsewhere.
In these studies the 1$d$ problem was transient.
More recently, coupled models discretized entirely by FEM were applied to
study cancer therapies, see e.g. \cite{nabil2015modelling} and references
therein. The mathematical foundations of these works are rooted in the seminal
contributions of \cite{d2008coupling}, \cite{d2012finite} where well-posedness
of the following problem is analyzed
\begin{equation}\label{kam_eq:daq}
\begin{aligned}
   -\Delta_{\Omega} u + (\hat{u}-\PI u)  \delta_{\Gamma} = f\delta_{\Gamma}   &\quad\text{ in } \Omega\\
   -\Delta_{\Gamma} \hat{u} - (\hat{u}-\PI u)  = g &\quad\text{ in } \Gamma.
\end{aligned}
\end{equation}

The presence of the measure term in \eqref{kam_eq:daq} requires the use of non-standard
spaces in the analysis. In \cite{d2008coupling}, the variational formulation of
the problem is proven to be well-posed using weighted Sobolev spaces. In particular, the
solution $u$ is sought in $H^1_{\alpha}(\Omega)$, $\alpha > 0$ while the test
functions $v$ for the first equation in \eqref{kam_eq:daq} are taken from $H^1_{-\alpha}(\Omega)$.
With this choice, the right-hand side $\langle f\delta_{\Gamma}, v \rangle$ as well as the trace operator
$v\mapsto v|_{\Gamma}$ are well defined, while the reduced regularity of $u$ is sufficient
to for the average $\PI u$ to make sense. As shown in \cite{d2012finite}, use of FEM
for the formulation in weighted spaces yields optimal rates if the computational
mesh is gradually refined towards $\Gamma$ (graded meshes).

Another approach to the analysis of \eqref{kam_eq:daq} has recently been suggested
in the numerical study \cite{cattaneo2015numerical}. Building on the
analysis of \cite{koppl2014optimal} for the elliptic problem with a 0 dimensional
Dirac right-hand side, the wellposedness of the problem was shown with
trial spaces $W^{1, p}(\Omega)$, $p=3-\tfrac{d}{2}$ and test spaces
$W^{1, q}(\Omega)$, $p^{-1}+q^{-1}=1$, and quasi-optimal error estimates for FEM
shown in the norms which excluded a fixed neighborhood of $\Gamma$ of radius $R$.

In studying AD or EPR, the physical parameters may vary across several orders
of magnitude while small or large time steps can be desirable depending on 
the time scales of interest. The solution algorithm for the employed model
equations is thus required to be robust with respect to these parameters 
as well independent of the discretization. For (the transient version of)
\eqref{kam_eq:daq} the construction of such algorithms is complicated by
the non-standard spaces on the domain $\Omega$.

A potential remedy for the
problem can be introduction of a Lagrange multiplier which enforces the
coupling between the domains with the goal of confining the non-standard
spaces to the smaller domain $\Gamma$. This idea has been used by \cite{kuchta2016preconditioners}
to analyze robust preconditioners for 2$d$-1$d$ coupled problems based on
operators in fractional Sobolev spaces, in particular, $(-\Delta_{\Gamma})^{\nhalf}$,
while numerical experiments reported in \cite{kuchta2016preconditioning}
suggest that for suitable exponents $(-\Delta_{\Gamma})^s$ defines a preconditioner for
the Schur complement of a 3$d$-1$d$ coupled system with a \emph{trace} constraint\footnote{
  Note that in \eqref{kam_eq:daq} and \eqref{kam_eq:strongform} the constraint/coupling 
  is defined in terms of a surface averaging operator $\PI$.
}.
We note that in both cases off-the-shelf methods were used as preconditioners for
the operators on $\Omega$.

The system \eqref{kam_eq:strongform} includes an additional variable $\lambda$
for the coupling constraint, cf. \eqref{kam_eq:daq}. We therefore aim 
to apply the techiques of \cite{kuchta2016preconditioning}, \cite{kuchta2016preconditioners} 
to construct a mesh-independent preconditioner for the problem, while 
the ideas of operator preconditioning \cite{mardal2011preconditioning} 
are used to ensure robustness with respect to the physical parameters and
the time-stepping.

The rest of the paper is organized as follows. Section \ref{kam_sec:prelim} 
identifies the structure of the preconditioner. In \S \ref{kam_sec:exp1} we discuss 
discretization of the proposed operator and report numerical experiments which 
demonstrate the robust properties. In \S \ref{kam_sec:exp2} the system 
\eqref{kam_eq:strongform} is used to model tissue perfusion using a realistic 
geometry of the rat cortex. Conclusions are finally summarized in \S \ref{kam_sec:concl}.


\section{Preconditioner for the coupled problem}\label{kam_sec:prelim}

In the following we let $\Omega$ be a bounded domain in $\reals^d$, $d=2$ or 3 and $\Gamma$ 
be a subdomain of $\Omega$ of dimension 1. By $L_2(D)$ we denote the space
of square-integrable functions over $D$ and $H^1(D)$ is the space of functions
with first order derivatives in $L^2(D)$. 


Discretizing \eqref{kam_eq:strongform} in time by backward-Euler discretization
the problem to be solved at each temporal level is of the form
\[
\A u = f
\] 
with
\begin{equation}\label{kam_eq:opA}
\A = \left[ 
\begin{array}{ccc} 
I - k D_\Omega \Delta_\Omega & 0 & k\beta \PI^*  \\   
0 & I - k D_\Gamma \Delta_\Gamma &  k\beta I  \\   
k\beta \PI & k\beta I  & -k.     
\end{array} 
\right]
\end{equation}
and $k$ being the time step size. Note that in order to obtain a symmetric
problem the operator $\mathcal{A}$ uses the adjoint $\PI^*$ of the \emph{averaging}
operator $\PI$ instead of the \emph{trace}, cf. \eqref{kam_eq:strongform}. 
The choice results in modeling error of order $\mathcal{O}(R)$. 

To motivate the structure of the preconditioner let us consider a $3\times 3$
matrix  
\[
A = \left[ 
\begin{array}{ccc} 
1 +\alpha_1 & 0 & \beta_1 \\   
0 & 1 + \alpha_2 & \beta_2  \\   
\beta_1  & \beta_2  & -\gamma
\end{array} 
\right]
\]
where $\alpha_1$, $\alpha_2$, $\beta_1$, $\beta_2$, and $\gamma$ are
assumed to be positive. It can then be shown that with
\[
B = \left[ 
\begin{array}{ccc} 
(1 +\alpha_1)^{-1} & 0 & 0 \\   
0 & (1 + \alpha_2)^{-1} & 0  \\   
0 & 0 & (\gamma + \beta_1^2 + \beta_2^2)^{-1} + (\gamma + \beta_1^2/\alpha_1 + \beta_2^2/\alpha_2)^{-1} 
\end{array} 
\right]
\]
the condition number of $B A$ is bounded independent of the parameters.
With this in mind we propose that $\mathcal{A}$ can be preconditioned by
a block-diagonal operator
\begin{equation}\label{kam_eq:opB}
\B = \left[ 
\begin{array}{ccc} 
(I - k D_\Omega \Delta_\Omega)^{-1} & 0 & 0 \\   
0 & (I - k D_\Gamma \Delta_\Gamma)^{-1} &  0  \\   
0 & 0 & S     
\end{array} 
\right], 
\end{equation}
where
\begin{equation}\label{kam_eq:opS}
  S=S^{-1}_1 + S_2^{-1}
\end{equation}
and
\begin{equation}\label{kam_eq:opS12}
  \begin{aligned}
    S_1 &= \gamma I + (k\beta)^2\PI \PI^*+ (k\beta)^2I,\\
    S_2  &=  \gamma I + (k \beta)^2 \PI(-k D_\Omega \Delta_\Omega)^{-1}\PI^* + (k \beta)^2(-k D_\Gamma \Delta_\Gamma)^{-1}.\\
    \end{aligned}
\end{equation}

The operator $\mathcal{B}$ could be rigorously derived within operator preconditioning %
\cite{mardal2011preconditioning} as a Riesz map preconditioner for $\mathcal{A}$
viewed as an isomorphism from $V(\Omega)\times \hat{V}(\Gamma) \times Q(\Gamma)$
to its dual space. In \cite{kuchta2016preconditioners} the framework
was applied to a system of two elliptic problems coupled by a 2$d$-1$d$
constraint. For \eqref{kam_eq:strongform} extension of the analysis to parabolic
problems would be required. In \cite{mardal2004uniform} robust preconditioners
for time-dependent Stokes problem were analyzed as operators between sums of
(parameter) weighted Sobolev spaces. Similarly, the structure of \eqref{kam_eq:opS}
suggests that $Q=Q_1 + Q_2$ with $Q_1$, $Q_2$ being suitable interpolation
spaces. However, here we shall not  justify $Q_1$ and $Q_2$
(and the preconditioner $\mathcal{B}$) theoretically. Instead, $Q_1$, $Q_2$
are characterized and robustness of $\mathcal{B}$ is demonstrated by numerical
experiments.

\section{Discrete preconditioner}\label{kam_sec:exp1}
Considering $\eqref{kam_eq:opB}$, both $(I - k D_\Omega \Delta_{\Omega})^{-1}$ and $(I - k D_\Gamma \Delta_\Gamma)^{-1}$  
can be realized with off-the-shelf multilevel algorithms and the crucial
question is thus how to construct $S$ efficiently. Note that assembling
$S_1$ and in particular $S_2$ might be too costly or even prohibitive, cf.
$(-\Delta_{\Omega})^{-1}$ in $S_2$. However, following \eqref{kam_eq:opS12} the
preconditioner can be realized if operators spectrally equivalent to $\PI (-\Delta_{\Omega})^{-1} \PI^*$
and $\PI \PI^*$ are known and if the inverse (action) of the resulting
approximations to $S_1$ and $S_2$ is inexpensive to compute.

\subsection{Auxiliary operators}

If $\Omega\subset\reals^2$ and $\PI$ is understood
as the trace operator, the mapping properties of trace as a bounded surjective
operator $H^1(\Omega)\rightarrow H^{\half}(\Gamma)$ can be used to show that
$\PI (-\Delta_{\Omega})^{-1} \PI^*$ is spectrally equivalent with $(-\Delta_{\Gamma})^{\nhalf}$.
At the same time $\PI \PI^*=\PI(-\Delta)^0_{\Omega}\PI^*$ requires characterizing
the space of traces of functions in $L^2{(\Omega)}$. We shall demonstrate
by a numerical experiment that a spectrally equivalent operator is here
provided by an operator $h^{-1} I_h$.

\begin{figure}[t]
  \begin{center}
  \includegraphics[width=0.35\textwidth, scale=0.6]{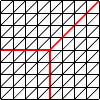}
    \hspace{10pt}
  \includegraphics[width=0.45\textwidth, clip]{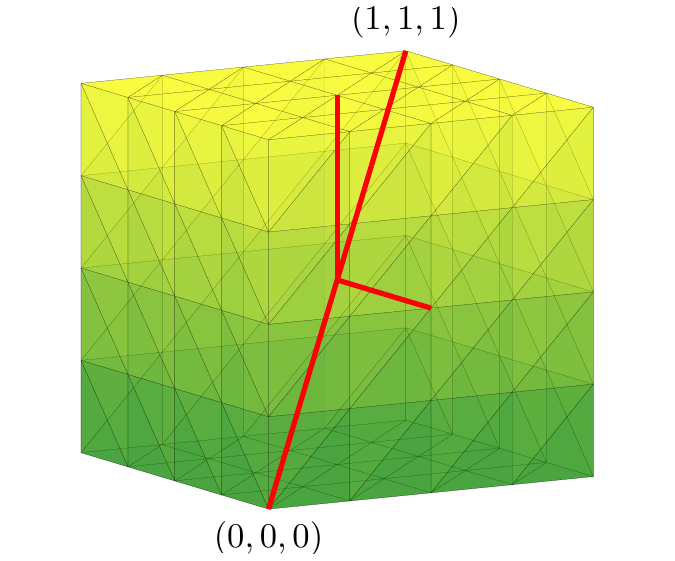}
\end{center}
\caption{Geometries used in preconditioning numerical experiments.
  The domain is $\Omega=\left[0, 1\right]^d$ while in order to prevent
  symmetries $\Gamma$ (pictured in red) always features a branching point.
  Triangulation of $\Gamma$ is made up of edges of the cells that triangulate
  $\Omega$.
}
\label{kam_fig:domains}
\end{figure}

Let $\Omega_h, \Gamma_h$ be triangulations of $\Omega$ and $\Gamma$ such that $\Gamma_h$ consists
of a subset of edges of the elements $\Omega_h$, cf. Figure \ref{kam_fig:domains}.
Further let $V_h, Q_h$ be finite element spaces of continuous linear Lagrange
elements on $\Omega_h$ and $\Gamma_h$ respectively. Finally we consider the
eigenvalue problem: Find $\lambda\in\reals$, $u\in V_h$, $p\in Q_h$ such
that
\begin{equation}\label{kam_eq:eigI}
  \begin{aligned}
    \int_{\Omega} u v + \int_{\Gamma} p \PI v &= \lambda  \int_{\Omega} u v & \quad \forall v\in V_h,\\ 
    \int_{\Gamma} q \PI u                   &= \lambda \int_{\Gamma}h^{-1} p q &\quad \forall q\in Q_h.
  \end{aligned}
\end{equation}
Table \ref{kam_tab:simple} shows the spectral condition number $\kappa=\max{\lvert \lambda \rvert}/\min{\lvert {\lambda} \rvert}$
of the linear systems \eqref{kam_eq:eigI}. For all the considered resolutions $h$ the
value of $\kappa$ is bounded. 

As the mapping properties of $\PI \PI^*$ and $\PI (-\Delta_{\Omega})^{-1} \PI^*$ in case
$\Omega\subset\reals^3$ are not trivially obtained from the continuous analysis,
we again resort to finding the suitable approximations by numerical experiments. Similar to the two dimensional
case, the first operator with $\PI$ having the constant radius $R=0.02$ is
found to be spectrally equivalent with $h^{-1} I_h$. Following \cite{kuchta2016preconditioning},
an approximation to $\PI (-\Delta_{\Omega})^{-1} \PI^*$ is searched for as a
suitable power $s < 0$ of $(-\Delta_{\Gamma} + I_{\Gamma})$. More precisely,
we look for the exponent yielding the most $h$-stable condition number of
the eigenvalue problem:
Find $\lambda\in\reals$, $u\in V_h$, $p\in Q_h$ such that
\begin{equation}\label{kam_eq:eigA}
  \begin{aligned}
    \int_{\Omega} \nabla u\cdot \nabla v + \int_{\Gamma} p \PI v &= \lambda  \int_{\Omega} \nabla u \nabla v + u v & \quad \forall v\in V_h,\\ 
    \int_{\Gamma} q \PI u                   &= \lambda \int_{\Gamma} p (-\Delta_{\Gamma} + I_{\Gamma})^s q &\quad \forall q\in Q_h.
  \end{aligned}
\end{equation}
In \eqref{kam_eq:eigA} the powers are computed using the spectral decomposition of the
operator $-\Delta_{\Gamma}+I_{\Gamma}$.

We shall not present here the results for the entire optimization problem
and only report on the optimum which is found to be $s=-0.55$. For this value
Table \ref{kam_tab:simple} shows the history of the condition numbers of system
\eqref{kam_eq:eigA} using again $R=0.02$. There is a slight growth by a constant
increment on the finer meshes, however, the final condition number is comparable
with that obtained on the coarsest mesh. We note that we have not investigated the
behaviour of the exponent on different curves or with variable radius. This subject
is left for future works.
%
\begin{table}[tb]
\small{
\begin{center}
  \begin{tabular}{c| >{\hspace{4pt}}c >{\hspace{4pt}}c >{\hspace{4pt}}c >{\hspace{4pt}}c >{\hspace{4pt}}c >{\hspace{4pt}}c ||
                     >{\hspace{4pt}}c >{\hspace{4pt}}c >{\hspace{4pt}}c >{\hspace{4pt}}c >{\hspace{4pt}}c >{\hspace{4pt}}c}
    \hline
    \multirow{3}{*}{} & \multicolumn{6}{c||}{$\Omega\subset\mathbb{R}^2$} & \multicolumn{6}{c}{$\Omega\subset\mathbb{R}^3$}\\
        \cline{2-13}
                                  & \multicolumn{6}{c||}{$1/h$} & \multicolumn{6}{c}{$1/h$}\\
    & 32 & 64 & 128 & 256 & 512 & 1024 & 4 & 8 & 16 & 32 & 64 & 128\\
    \hline
eq. \eqref{kam_eq:eigI} & 4.36 & 4.33 & 4.35 & 4.36 & 4.36 & 4.35 & 4.26 & 5.22 & 5.58 & 4.44 & 4.85 & 4.85\\
eq. \eqref{kam_eq:eigA} & 8.80 & 8.84 & 8.84 & 8.86 & 8.88 & 8.87 & 11.04 & 8.25 & 7.44 & 9.02 & 10.30 & 11.25\\
    \hline
  \end{tabular}
\end{center}
}
\caption{Spectral condition numbers of the eigenvalue problems related to
  approximations of $\PI \PI^*$(eq \eqref{kam_eq:eigI}) and $\PI(-\Delta_{\Omega})\PI^*$(eq \eqref{kam_eq:eigA}).
  In the two-dimensional case \eqref{kam_eq:eigA} uses $s=-\tfrac{1}{2}$ in agreement
  with the mapping properties of the continuous trace operator. Results for $s=-0.55$ are reported
  in the three dimensional case. On the finest triangulation $\dim V_h \sim 10^6$
  and $\dim Q_h \sim 10^3$ when $d=2$ and $\dim Q_h \sim 10^2$ for $d=3$. 
}
\label{kam_tab:simple}
\end{table}


\subsection{Discrete preconditioner for the coupled problem}

Applying the proposed preconditioner \eqref{kam_eq:opB} of the coupled
3$d$-1$d$ problem \eqref{kam_eq:opA} requires evaluating the inverse of
operators $S_1$ and $S_2$ in \eqref{kam_eq:opS12}. The former is readily
computed since, due to the suggested equivalence of $\PI \PI^*$ and $h^{-1} I_h$,
the matrix representation of $S_1$ is essentialy a rescaled mass matrix.
For $S_2$ we show that if $(-\Delta_{\Gamma} + I_{\Gamma})^s$ is computed from the spectral
decomposition then the inverse $S^{-1}_{2}$ can be computed in a closed form.

Let $A$, $M$ be the $n\times n$ matrix representations
of Galerkin approximations of $-\Delta_{\Gamma}+I_{\Gamma}$ and $I$ in the space $Q_h\subset H^1(\Gamma)$.
Following \cite{kuchta2016preconditioners} the matrix representation of
$(-\Delta_{\Gamma} + I_{\Gamma})_h^s$ is $H_s = M U \Lambda ^s (MU)'$ where the matrices $U$, $\Lambda$
solve the generalized eigenvalue problem $AU=MU\Lambda$ such that $U'MU=I$.  
Using $H_s$ it is easily established that the matrix representation of $S_2$ is
\[
  \gamma H_0 + (k\beta)^2(k D_{\Omega})^{-1}H_{-\frac{1}{2}} + (k\beta)^2(k D_{\Gamma})^{-1}H_{-1}.
\]
As $H^{-1}_s=U\Lambda^{-s}U'$ the inverse of the $S_2$ matrix is given by
\begin{equation}\label{kam_eq:S2inv}
U\left(\gamma \Lambda^0 + (k\beta)^2(k D_{\Omega})^{-1}\Lambda^{-\frac{1}{2}} + (k\beta)^2(k D_{\Gamma})^{-1}\Lambda^{-1}\right)^{-1}U'.
\end{equation}
%

Using the spectral decomposition, the cost of setting up the preconditioner
$S$ is determined by the cost of solving the generalized eigenvalue problem
for $U$ and $\Lambda$. This practically limits the construction to systems where
$\dim Q_h \sim 10^3$. However, for the problems considered further, this limitation
does not present an issue. In particular, the preconditioner can be setup on
the discretization of vasculature of the cortex tissue used in \S\ref{kam_sec:exp2}
which contains approximately twenty thousand vertices.

%
%

With a discrete approximation of $S$ we can finally address the question of
$\mathcal{B}$ being a robust preconditioner for \eqref{kam_eq:opA}. Motivated by 
\eqref{kam_eq:strongform} we do not vary all the paramaters, and instead, set $\gamma=1$
in $\mathcal{A}$. Morever, the conductivity on $\Omega$ is taken as unity and
only variable $D_{\Gamma} > 1$ is considered mimicking the expected faster
propagation along the one-dimensional domain. Finally, the time step $k$
and the coupling constant $\beta$ shall take values between $10^{-4}$ and $10^{-8}$.

For a fixed choice of parameters $D_{\Gamma}$, $\beta$, $k$, the preconditioned problem
$\mathcal{B}\mathcal{A}x = \mathcal{B}f$ is considered on geometries from
Figure \ref{kam_fig:domains} and discretized with continuous linear Lagrange
elements. The resulting linear system is then solved with the MinRes method
where the iterations are stopped once the preconditioned residual norm is
less than $10^{-10}$ in magnitude.

The observed iteration counts are reported in Table \ref{kam_tab:iters}.
For both 2$d$-1$d$ and 3$d$-1$d$ coupled problems the iterations can be seen
to be bounded in the discretization parameter. Moreover the preconditioner
performs almost uniformly in the considered ranges of $D_{\Gamma}$ and $\beta$
while there is a clear boundedness in $k$ as well. We note that these conclusions
are not significantly altered if the ranges for $k$ and $\beta$ are extended to 1.

\begin{table}
\small{
\begin{center}
  \begin{tabular}{c|c|c|| >{\hspace{6pt}}c >{\hspace{6pt}}c >{\hspace{6pt}}c >{\hspace{6pt}}c >{\hspace{6pt}}c >{\hspace{6pt}}c
                      ||| >{\hspace{6pt}}c >{\hspace{6pt}}c >{\hspace{6pt}}c >{\hspace{6pt}}c >{\hspace{6pt}}c >{\hspace{6pt}}c}
    \hline
    \multirow{3}{*}{$D_{\Gamma}$} & \multirow{3}{*}{$\beta$} & \multirow{3}{*}{$k$} & \multicolumn{6}{c|||}{2$d$-1$d$} & \multicolumn{6}{c}{3$d$-1$d$}\\
        \cline{4-15}
                                  & & & \multicolumn{6}{c|||}{$1/h$} & \multicolumn{6}{c}{$1/h$}\\
    & & &  32 & 64 & 128 & 256 & 512 & 1024 & 4 & 8 & 16 & 32 & 64 & 128\\
    \hline
\multirow{9}{*}{$10^{0}$} & \multirow{3}{*}{$10^{-8}$} & $10^{-8}$ & 11 & 10 & 10 & 8 & 7 & 7 & 8 & 13 & 12 & 11 & 11 & 10\\
  &   & $10^{-6}$ & 15 & 13 & 12 & 10 & 8 & 8 & 10 & 14 & 16 & 16 & 15 & 13\\
  &   & $10^{-4}$ & 16 & 11 & 12 & 13 & 13 & 14 & 12 & 16 & 18 & 16 & 13 & 14\\
\cline{3-15}
  & \multirow{3}{*}{$10^{-6}$} & $10^{-8}$ & 11 & 10 & 10 & 8 & 7 & 7 & 8 & 13 & 12 & 11 & 11 & 10\\
  &   & $10^{-6}$ & 15 & 13 & 12 & 10 & 8 & 8 & 10 & 14 & 16 & 16 & 15 & 13\\
  &   & $10^{-4}$ & 15 & 11 & 12 & 13 & 13 & 14 & 12 & 16 & 18 & 16 & 13 & 14\\
\cline{3-15}
  & \multirow{3}{*}{$10^{-4}$} & $10^{-8}$ & 11 & 10 & 10 & 8 & 7 & 7 & 8 & 13 & 12 & 11 & 11 & 10\\
  &   & $10^{-6}$ & 15 & 13 & 12 & 10 & 8 & 8 & 10 & 14 & 16 & 16 & 15 & 13\\
  &   & $10^{-4}$ & 16 & 11 & 12 & 13 & 13 & 14 & 12 & 16 & 18 & 16 & 13 & 14\\
\cline{2-15}
\multirow{9}{*}{$10^{2}$} & \multirow{3}{*}{$10^{-8}$} & $10^{-8}$ & 11 & 10 & 10 & 9 & 7 & 7 & 8 & 13 & 12 & 11 & 11 & 10\\
  &   & $10^{-6}$ & 15 & 13 & 12 & 10 & 8 & 8 & 10 & 14 & 16 & 16 & 15 & 13\\
  &   & $10^{-4}$ & 15 & 11 & 12 & 13 & 13 & 13 & 11 & 16 & 18 & 13 & 13 & 14\\
\cline{3-15}
  & \multirow{3}{*}{$10^{-6}$} & $10^{-8}$ & 11 & 10 & 10 & 9 & 7 & 7 & 8 & 13 & 12 & 11 & 11 & 10\\
  &   & $10^{-6}$ & 15 & 13 & 12 & 10 & 8 & 8 & 10 & 14 & 16 & 16 & 15 & 13\\
  &   & $10^{-4}$ & 15 & 11 & 12 & 13 & 13 & 13 & 11 & 16 & 18 & 13 & 13 & 14\\
\cline{3-15}
  & \multirow{3}{*}{$10^{-4}$} & $10^{-8}$ & 11 & 10 & 10 & 9 & 7 & 7 & 8 & 13 & 12 & 11 & 11 & 10\\
  &   & $10^{-6}$ & 15 & 13 & 12 & 10 & 8 & 8 & 10 & 14 & 16 & 16 & 15 & 13\\
  &   & $10^{-4}$ & 15 & 11 & 12 & 13 & 13 & 13 & 12 & 16 & 18 & 13 & 13 & 14\\
\cline{2-15}
\multirow{9}{*}{$10^{4}$} & \multirow{3}{*}{$10^{-8}$} & $10^{-8}$ & 11 & 10 & 10 & 8 & 7 & 7 & 8 & 13 & 12 & 11 & 11 & 10\\
  &   & $10^{-6}$ & 15 & 13 & 12 & 10 & 8 & 8 & 10 & 14 & 16 & 16 & 15 & 13\\
  &   & $10^{-4}$ & 15 & 11 & 12 & 13 & 12 & 10 & 11 & 16 & 18 & 13 & 13 & 14\\
\cline{3-15}
  & \multirow{3}{*}{$10^{-6}$} & $10^{-8}$ & 11 & 10 & 10 & 8 & 7 & 7 & 8 & 13 & 12 & 11 & 11 & 10\\
  &   & $10^{-6}$ & 15 & 13 & 12 & 10 & 8 & 8 & 10 & 14 & 16 & 16 & 15 & 13\\
  &   & $10^{-4}$ & 15 & 11 & 12 & 13 & 12 & 10 & 11 & 16 & 18 & 13 & 13 & 14\\
\cline{3-15}
  & \multirow{3}{*}{$10^{-4}$} & $10^{-8}$ & 11 & 10 & 10 & 8 & 7 & 7 & 8 & 13 & 12 & 11 & 11 & 10\\
  &   & $10^{-6}$ & 15 & 13 & 12 & 10 & 8 & 8 & 10 & 14 & 16 & 16 & 15 & 13\\
  &   & $10^{-4}$ & 15 & 11 & 12 & 13 & 12 & 10 & 11 & 16 & 18 & 13 & 13 & 14\\
\cline{2-15}
\multirow{9}{*}{$10^{6}$} & \multirow{3}{*}{$10^{-8}$} & $10^{-8}$ & 11 & 10 & 10 & 8 & 7 & 7 & 8 & 13 & 12 & 11 & 11 & 10\\
  &   & $10^{-6}$ & 15 & 13 & 12 & 10 & 8 & 8 & 10 & 14 & 16 & 16 & 15 & 13\\
  &   & $10^{-4}$ & 14 & 10 & 10 & 9 & 8 & 5 & 11 & 15 & 18 & 13 & 11 & 11\\
\cline{3-15}
  & \multirow{3}{*}{$10^{-6}$} & $10^{-8}$ & 11 & 10 & 10 & 8 & 7 & 7 & 8 & 12 & 12 & 11 & 11 & 10\\
  &   & $10^{-6}$ & 15 & 13 & 12 & 10 & 8 & 8 & 10 & 14 & 16 & 16 & 15 & 13\\
  &   & $10^{-4}$ & 14 & 10 & 10 & 9 & 8 & 5 & 11 & 15 & 18 & 13 & 10 & 11\\
\cline{3-15}
  & \multirow{3}{*}{$10^{-4}$} & $10^{-8}$ & 11 & 10 & 10 & 8 & 7 & 7 & 8 & 12 & 12 & 11 & 11 & 10\\
  &   & $10^{-6}$ & 15 & 13 & 12 & 10 & 8 & 8 & 10 & 14 & 16 & 16 & 15 & 13\\
  &   & $10^{-4}$ & 14 & 10 & 10 & 9 & 8 & 5 & 11 & 16 & 18 & 13 & 10 & 11\\
\cline{3-15}
\hline
\end{tabular}
\end{center}
\caption{Number of iterations of MinRes method on \eqref{kam_eq:opA} using  \eqref{kam_eq:opB} as
  preconditioner with $S$ approximated using \eqref{kam_eq:S2inv}.
  (Left) 2$d$-1$d$ coupled problem and (right) 3$d$-1$d$ coupled problem
  from Figure \ref{kam_fig:domains} are considered.}
\label{kam_tab:iters}
}
\end{table}

Having demonstrated the numerical stability of our model, we next test it
on the same problem as \cite{tofts2012modelling},
namely a bloodborne tracer perfusing and later being cleared from tissue.
While \cite{tofts2012modelling} considers this problem on a macroscopic scale,
we model it on the micro-scale, where individual blood vessels can be resolved
as part of our 1$d$ domain.

\section{Perfusion experiment}\label{kam_sec:exp2}
In \cite{sakadvzic2014large}, a $0.7 \text{mm} \times 0.7 \text{mm} \times 0.7 \text{mm}$
piece of mouse brain microvasculature was imaged using two-photon microscopy.
To obtain a realistic geometry for our model, we used this data to generate a
3$d$ mesh of the extravascular space in which vessel segments corresponded
to 1$d$ mesh edges. The radius of the blood vessels is used as the radius $R$ 
in the definition of the averaging operator $\PI$, and ranged between $1$ and $15$ micron.

To model a small region of tissue being perfused by a bloodborne tracer,
we use initial conditions of $u, \hat{u} = 0$,
and a boundary condition of $\hat{u} = 1 \: \frac {\text{mol}} {\text{L}}$ on the part of the
boundary corresponding to inlet vessels. To model clearance, the inlet
boundary condition was swapped from $\hat{u} = 1$ to $\hat{u} = 0$ after
a third of the simulation time had passed.

As parameters, we use $D_{\Gamma} = 6.926\times 10^7 \: \micron^2 / \text{s}$, $D_{\Omega} =  1.87\times 10^2 \: \micron^2 / \text{s}$
and $\beta = 50 \: \micron / \text{s}$. We used $k = 1 \text{s}$ after
verifying that reducing it did not significantly affect our results.
In this experiment, it was unneccessary to use the preconditioner
described in \S \ref{kam_sec:exp1} since the problem size was small
enough to allow for use of a direct solver.

\begin{figure}
  \begin{center}
    \includegraphics[width=0.4\textwidth]{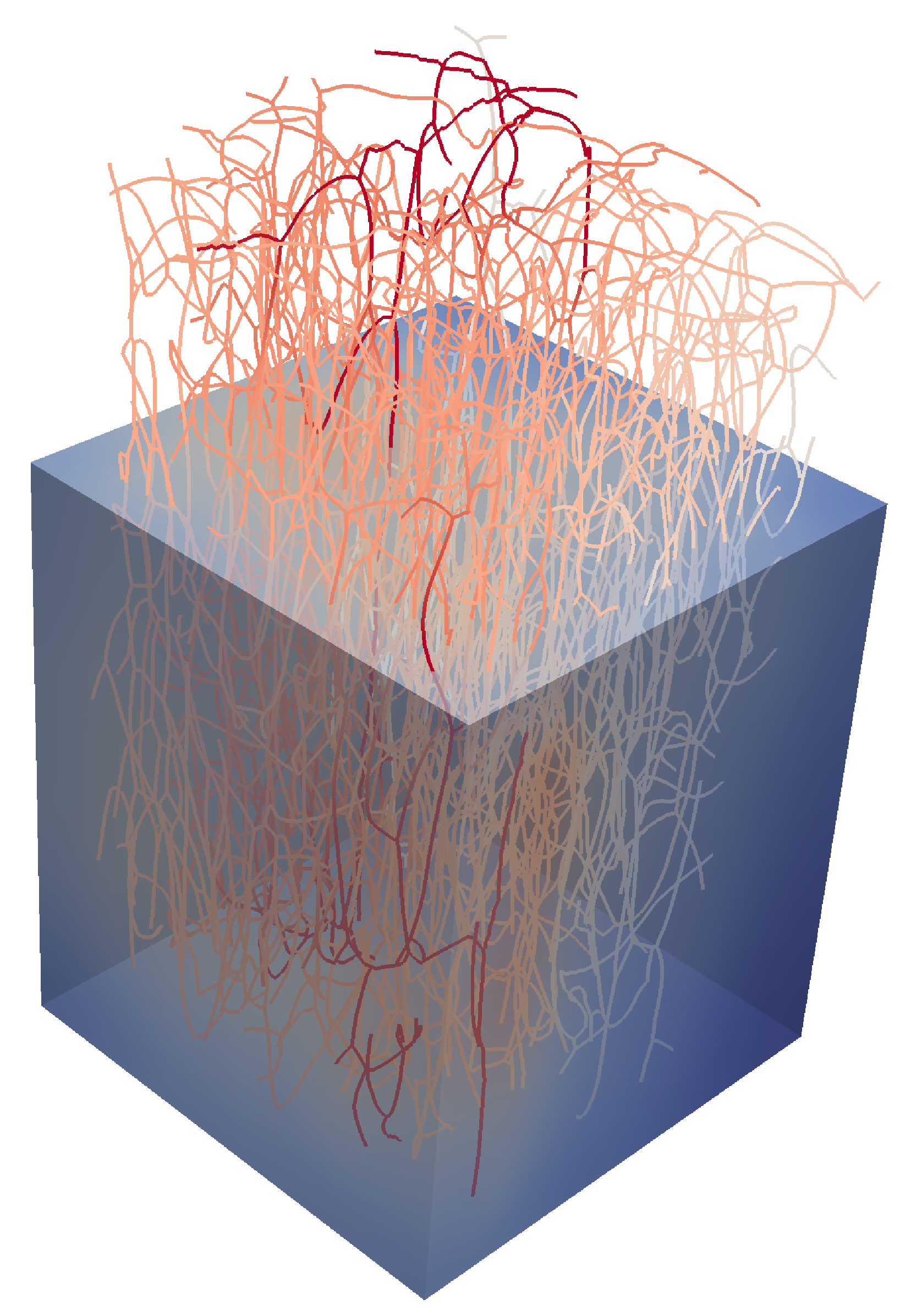}
    \includegraphics[width=0.4\textwidth]{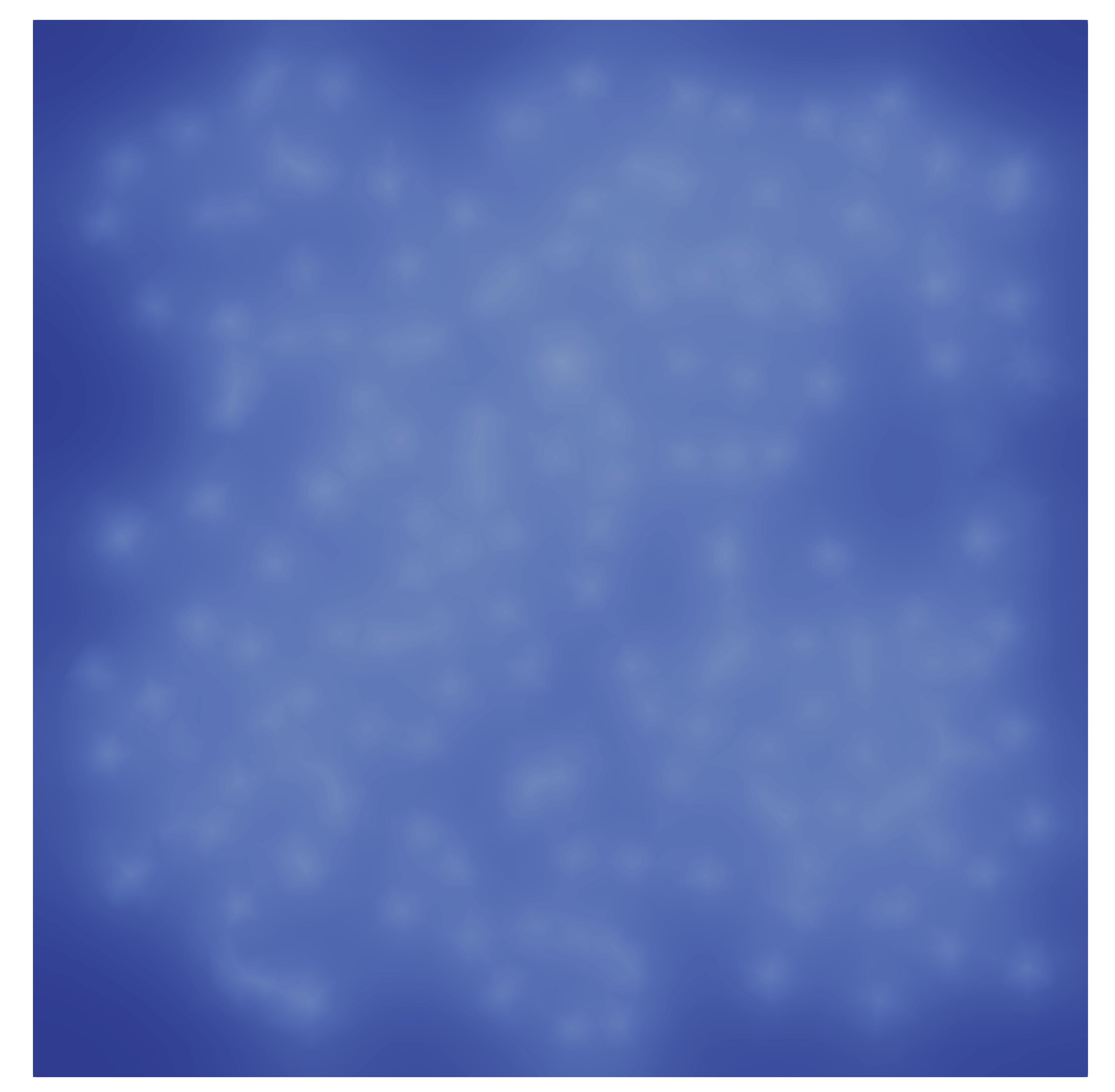}
    \caption{Example results shown on a (left) clip of the 3$d$ domain and (right) on a slice.
      Notice the 'halos' of increased concentration immediately around the vessels.}
    \label{kam_fig:domain}
  \end{center}
\end{figure}

Tofts \cite{tofts2012modelling} assumes a relation
\begin{align}
  \label{kam_eq:D}
  \frac{\partial C_t} {\partial t} = \frac {K_{\text{trans}}} {\nu} (C_v - C_t) 
\end{align}
between the pixel tissue concentration $C_t$ and the pixel vessel concentration $C_v$
for some constant $K_{\text{trans}}$. Here, $\nu$ is the vascular volume fraction, which in our geometry is about  $0.76 \%$. Our geometry is of
a size comparable to a single pixel in \cite{tofts2012modelling}, so
$C_t$ and $C_v$ correspond to the normalized averages
\[
  C_t = \frac{\int_{\Omega} u} {\int_{\Omega} 1}
  \quad\mbox{and}\quad
  C_v = \frac{\int_{\Gamma} \pi R^2 \hat{u}} {\int_{\Gamma} \pi R^2}.
\]

We computed $K_{\text{trans}}$ by solving for $u, \hat{u}$ using our model, and then computing
$C_t, C_v$ as given above, and defining $K_{\text{trans}}$ such that equation \eqref{kam_eq:D}
holds at each time point. This makes $K_{\text{trans}}$ a function of time, with units $\text{seconds}^{-1}$.

\begin{figure}
  \begin{center}                
    \includegraphics[width=0.49\textwidth]{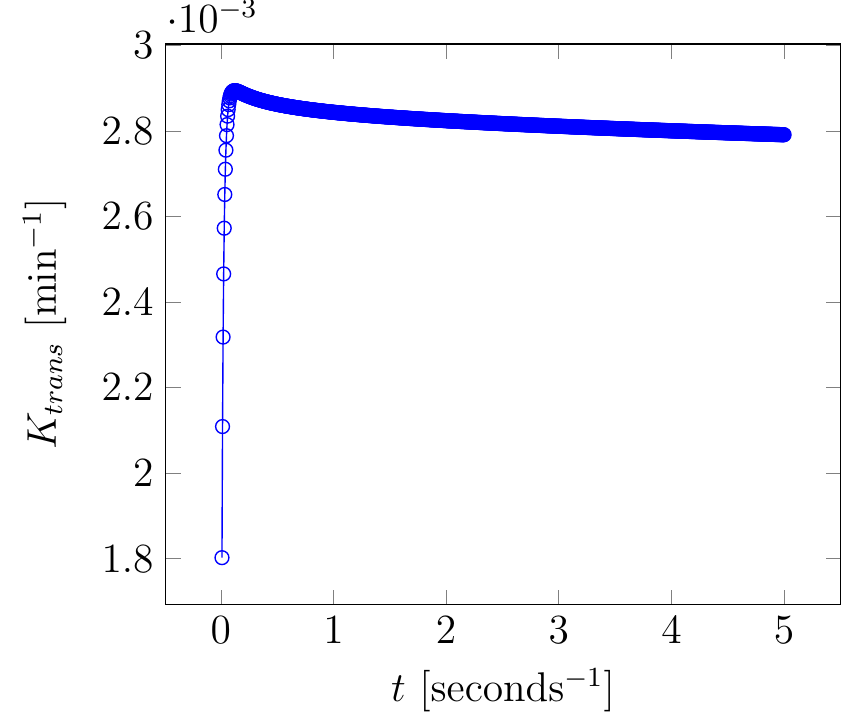}
    \includegraphics[width=0.49\textwidth]{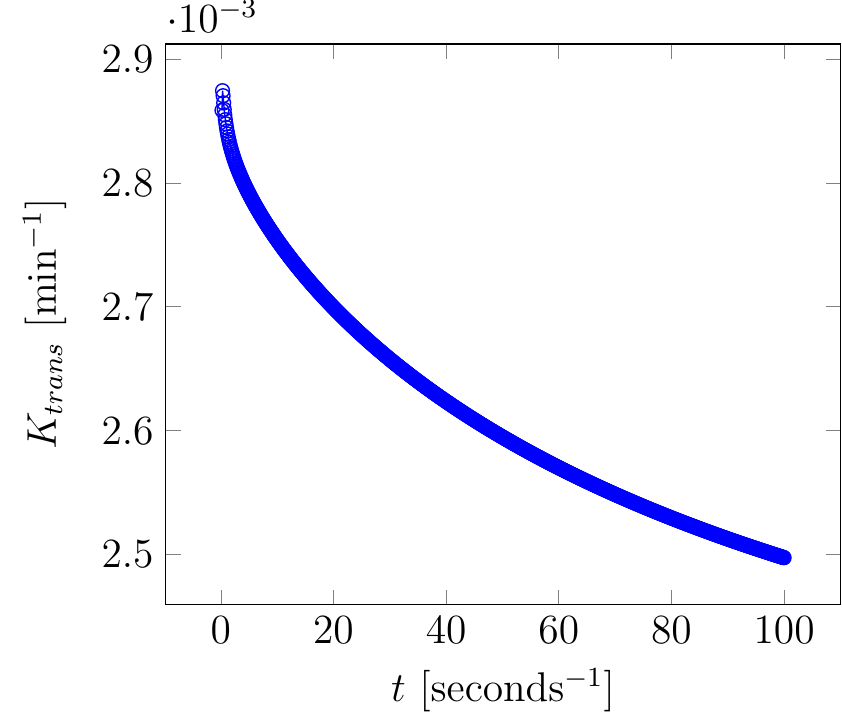}
    \caption{Behavior of $K_{\text{trans}}$ on short and long time scales.}
    \label{kam_fig:K_graph}
  \end{center}
\end{figure}
\subsection{Discussion of perfusion experiment}
Our value for $K_{\text{trans}}$ is not entirely constant. There is a small variation in time,
as perfusion seems to be somewhat faster when the extravascular space is
completely empty of the tracer. As it starts getting saturated, perfusion
slows down somewhat. This translates into $K_{\text{trans}}$ decreasing
by about $20\%$ over a period of about 5 minutes, from 0.0028 $\text{min}^{-1}$ to 0.0023 $\text{min}^{-1}$.

In \cite{zhang2012correlation}, $K_{\text{trans}}$ was estimated in healthy
and cancerous human brain tissue from MRI scans. In healthy tissue, they
estimate $K_{\text{trans}}$ to be between 0.003 $\text{min}^{-1}$ and 0.005 $\text{min}^{-1}$,
that is, slightly higher than our results. There are several possible
explanations for this difference. One might be that in our model, vascular transport
is modeled as exceptionally fast diffusion for convenience, whereas in
reality it occurs by convection. However,
in both cases 1$d$ transport is very fast compared to the 3$d$ transport and
the 1$d$-3$d$ exchange. Further, $K_{\text{trans}}$ is defined in terms of the 1$d$-3$d$ exchange
alone, so non-extreme variations in the 1$d$ transport seem unlikely to be relevant.

Another possibility might be that the data of \cite{zhang2012correlation} are taken from human brain tissue, while our vasculature is taken from a mouse brain tissue, likely from a different region of the brain. A third
reason might be our diffusion constants not exactly matching the tracer used
by \cite{zhang2012correlation}. 

In further work, it would be interesting to incorporate convective transport
into the model and see if better agreement with the experimental data is
observed. A suitable starting point here is \cite{cattaneo2014computational},
who derive a convection-diffusion type system (equations (3a), (3b))
by assuming that the blood flow $\hat{q}$ in a segment is laminar
and follows Poiseuille's law $$\hat{q} =  R^4 C \nabla {\hat{p}}.$$

\subsection{Parameter sensitivity analysis}
We carried out a rudimentary parameter sensitivity analysis by varying the three parameters $D_{\Gamma}, D_{\Omega}, \beta$ by a factor 10 and seeing how that affected the tissue concentration $C_t$. Specifically, we started from a baseline of \\ $D_{\Gamma}= 6.926 \cdot 10^8 \: \frac {\mu m^2} {s}, D_{\Omega}= 1.87 \cdot 10^2 \: \frac {\mu m^2} {s}, \beta = 5 \cdot 10^1 \: \text{s}^{\nhalf}$, and for each parameter, increased or decreased it by a factor 10. 

\begin{figure}
  \begin{center}
    \includegraphics[width=0.31\textwidth]{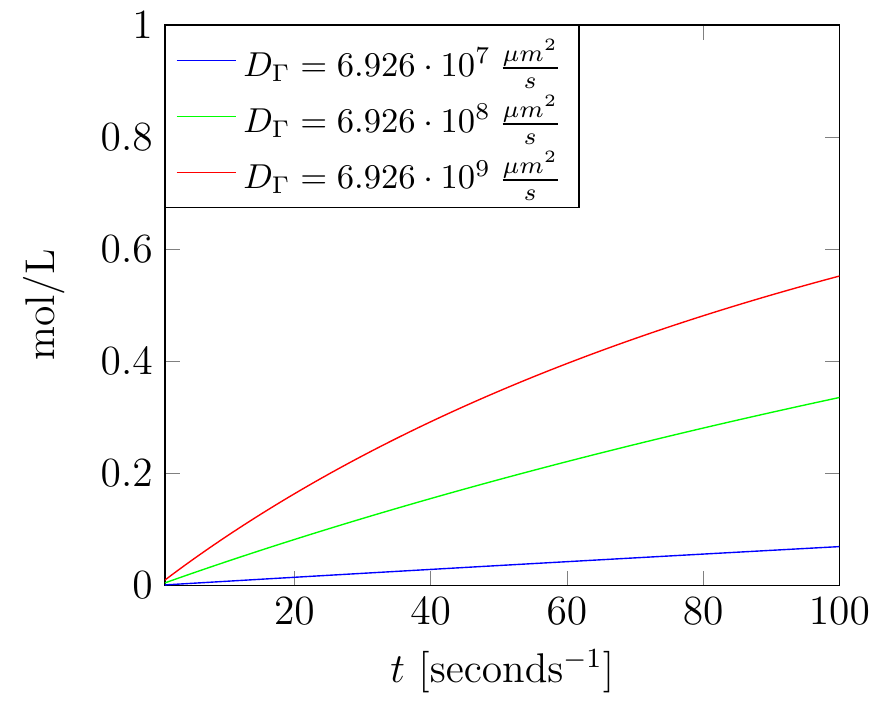}
    \includegraphics[width=0.31\textwidth]{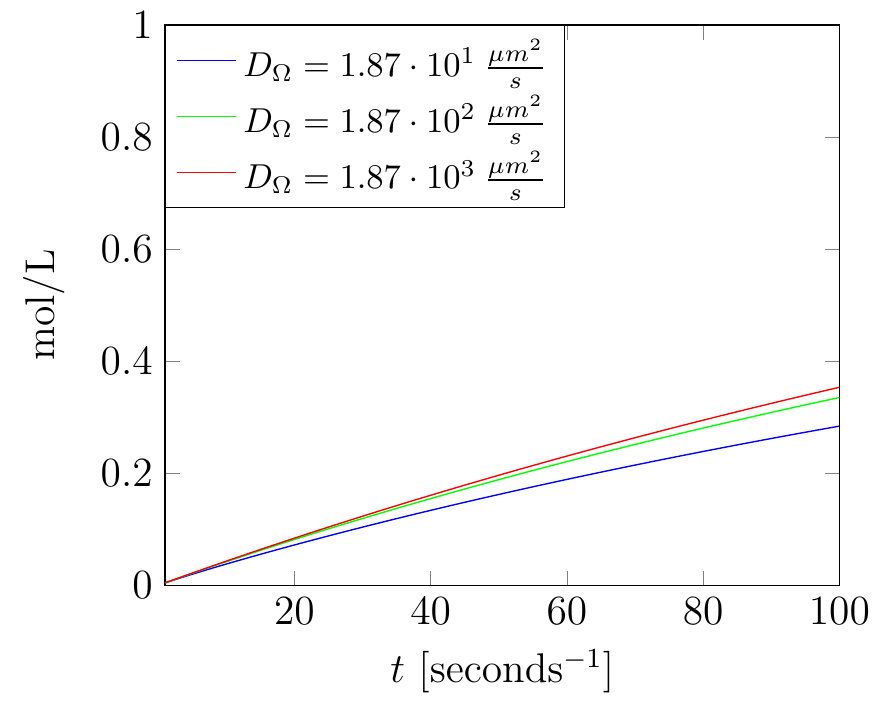}
    \includegraphics[width=0.31\textwidth]{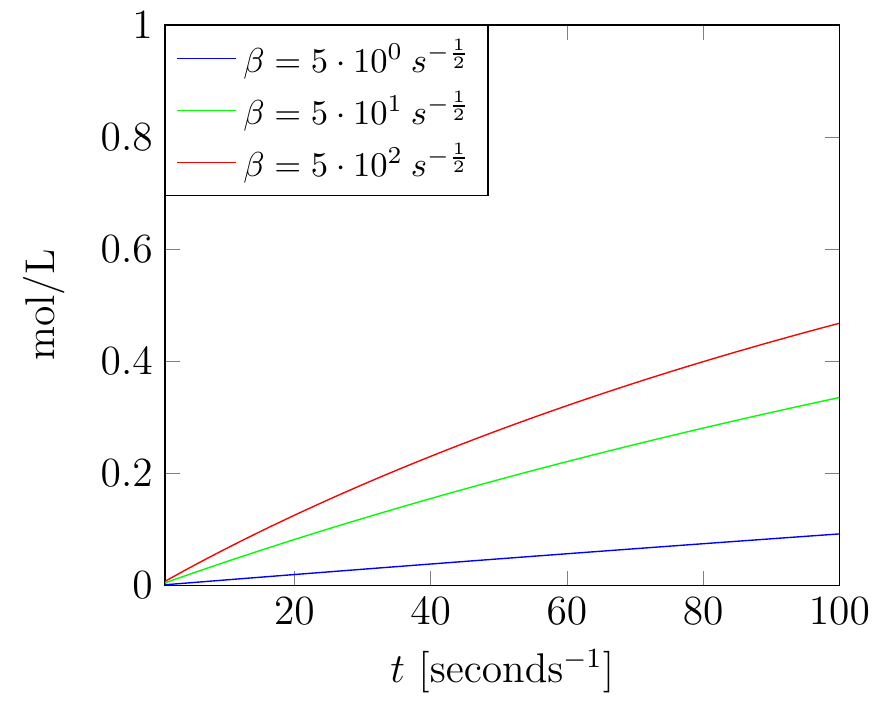}

    \caption{Plots of variation in $C_t$ when different parameters are varied. }
    \label{kam_fig:param_variation}
  \end{center}
\end{figure}

The results are shown in Figure \ref{kam_fig:param_variation}. They indicate that $C_t$ depends more strongly on $\beta$ and $D_{\Gamma}$ than on $D_{\Omega}$ for the set of parameters considered here. 

\section{Conclusions}\label{kam_sec:concl}
A coupled 3$d$-1$d$ system with an additional unknown enforcing the coupling
between the domains was used as a model of tissue perfusion. For the system
we proposed a robust preconditioner and demonstrated its properties through
numerical experiments. Further, we have shown that the model can be applied to a
physiological problem with reasonable results.


%
%
%
\bibliographystyle{vmams}  
\bibliography{refs}        

\end{document}